# New rapidly converging series representations for values of the Riemann zeta function $\varsigma(2n+1)$ and the Dirichlet beta function $\beta(2n)$

Donal F. Connon

dconnon@btopenworld.com

24 March 2010

**Abstract**

In this paper we derive rapidly converging series for Catalan's constant $G = \beta(2)$ and for Apéry's constant $\varsigma(3)$. These are set out below:

$$G = 3(1-\log 2) - \sum_{n=1}^{\infty} \frac{1}{n} \frac{1}{(16n^2-1)^2}$$

$$7\varsigma(3) = 12(5 - 4\log 2) - 20G + 16\sum_{n=1}^{\infty} \frac{1}{n} \frac{1}{(16n^2-1)^3}$$

The method may be easily generalised to produce new series representations for values of the Riemann zeta function $\varsigma(2n+1)$ and the Dirichlet beta function $\beta(2n)$.

**Catalan's constant**

We have the Weierstrass canonical form of the gamma function [13, p.1]

(1) $$\frac{1}{\Gamma(x)} = xe^{\gamma x} \prod_{n=1}^{\infty} \left\{ \left(1 + \frac{x}{n}\right) e^{-\frac{x}{n}} \right\}$$

and taking logarithms results in

(2) $$\log \Gamma(x) = -\log x - \gamma x - \sum_{n=1}^{\infty} \left[ \log\left(1 + \frac{x}{n}\right) - \frac{x}{n} \right]$$

The following well known identity for the digamma function $\psi(x)$ is easily derived by differentiating equation (2)

(3) $$\psi(x) = \frac{\Gamma'(x)}{\Gamma(x)} = -\frac{1}{x} - \gamma - \sum_{n=1}^{\infty} \left( \frac{1}{x+n} - \frac{1}{n} \right)$$

Since [13, p.14]

(4) $$\psi(1+x) = \psi(x) + \frac{1}{x}$$

we obtain

$$\psi(1+x) + \gamma = -\sum_{n=1}^{\infty}\left(\frac{1}{x+n} - \frac{1}{n}\right)$$

and we easily see that

$$\psi(1+x) + \psi(1-x) + 2\gamma = -\sum_{n=1}^{\infty}\left(\frac{1}{n+x} + \frac{1}{n-x} - \frac{2}{n}\right)$$

or equivalently

(5) $$\psi(1+x) + \psi(1-x) + 2\gamma = -2x^2 \sum_{n=1}^{\infty} \frac{1}{n} \frac{1}{n^2 - x^2}$$

which concurs with Prudnikov et al. [11, 5.1.25-13].

With $x = p/2$ we have

(6) $$\psi\left(1+\frac{p}{2}\right) + \psi\left(1-\frac{p}{2}\right) + 2\gamma = -2p^2 \sum_{n=1}^{\infty} \frac{1}{n} \frac{1}{4n^2 - p^2}$$

Ramanujan [21, p.29] showed that

(7) $$\sum_{n=1}^{\infty} \frac{1}{n} \frac{1}{16n^2 - 1} = 3\log 2 - 2$$

In fact, as reported in [3], this summation was the **first** problem submitted by Ramanujan in 1911 to the Journal of the Indian Mathematical Society. Ramanujan's result follows by letting $p = 1/2$ in (6).

We have [13, p.20]

$$\psi\left(\frac{3}{4}\right) = -\gamma + \frac{1}{2}\pi - 3\log 2 \qquad \psi\left(\frac{5}{4}\right) = -\gamma - \frac{1}{2}\pi - 3\log 2 + 4$$

and (7) consequently follows.

Ramanujan [2, p.26] also showed that

(8) $$\sum_{n=1}^{\infty} \frac{1}{n} \frac{1}{4n^2-1} = 2\log 2 - 1$$

This series is contained in [7, p.9,] where it is attributed to Bromwich's text [5, p.51]. A proof is also given in Ramanujan's Notebooks [2, p.26]. In [5] the proof is based on the definition of Euler's constant $\gamma$ and it is therefore quite apt that the series should be connected with the digamma function. The identity (8) may be obtained by letting $p=1$ in (6).

Differentiating (6) results in

(9) $$\frac{1}{2}\psi'\left(1+\frac{p}{2}\right) - \frac{1}{2}\psi'\left(1-\frac{p}{2}\right) = -4p^3 \sum_{n=1}^{\infty} \frac{1}{n} \frac{1}{(4n^2-p^2)^2} - 4p \sum_{n=1}^{\infty} \frac{1}{n} \frac{1}{4n^2-p^2}$$

and with $p=1$ we obtain

$$\frac{1}{2}\psi'\left(1+\frac{1}{2}\right) - \frac{1}{2}\psi'\left(\frac{1}{2}\right) = -4 \sum_{n=1}^{\infty} \frac{1}{n} \frac{1}{(4n^2-1)^2} - 4 \sum_{n=1}^{\infty} \frac{1}{n} \frac{1}{4n^2-1}$$

Since $\psi(1+x) - \psi(x) = \frac{1}{x}$ we have

(10) $$\psi'(1+x) - \psi'(x) = -\frac{1}{x^2}$$

and thus

$$\frac{1}{2}\psi'\left(1+\frac{1}{2}\right) - \frac{1}{2}\psi'\left(\frac{1}{2}\right) = -2$$

Hence using (8) we obtain

(11) $$\sum_{n=1}^{\infty} \frac{1}{n} \frac{1}{(4n^2-1)^2} = \frac{3}{2} - 2\log 2$$

which is reported in [7, p.9].

Kölbig [10] has shown that

(12) $$\psi^{(2k-1)}\left(\frac{1}{4}\right) - \psi^{(2k-1)}\left(\frac{3}{4}\right) = (2k-1)! 2^{4k} \beta(2k)$$

(13) $$\psi^{(2k)}\left(\frac{1}{4}\right)+\psi^{(2k)}\left(\frac{3}{4}\right)=-(2k)!2^{2k+1}(2^{2k+1}-1)\varsigma(2k+1)$$

where $\beta(s)$ is the Dirichlet beta function defined by

$$\beta(s)=\sum_{j=0}^{\infty}\frac{(-1)^j}{(2j+1)^s}$$

and $G=\beta(2)$ is called Catalan's constant.

In particular we have

(14) $$\psi^{(1)}\left(\frac{1}{4}\right)-\psi^{(1)}\left(\frac{3}{4}\right)=16G$$

and

(15) $$\psi^{(2)}\left(\frac{1}{4}\right)+\psi^{(2)}\left(\frac{3}{4}\right)=-112\varsigma(3)$$

With $p=1/2$ in (9) we get

$$\psi'\left(1+\frac{1}{4}\right)-\psi'\left(\frac{3}{4}\right)=-16\sum_{n=1}^{\infty}\frac{1}{n}\frac{1}{(16n^2-1)^2}-16\sum_{n=1}^{\infty}\frac{1}{n}\frac{1}{16n^2-1}$$

$$=-16\sum_{n=1}^{\infty}\frac{1}{n}\frac{1}{(16n^2-1)^2}-48\log 2+32$$

We have from (10)

$$\psi'\left(1+\frac{1}{4}\right)=\psi'\left(\frac{1}{4}\right)-16$$

and therefore we obtain

(16) $$G=3(1-\log 2)-\sum_{n=1}^{\infty}\frac{1}{n}\frac{1}{(16n^2-1)^2}$$

This series representation appears to be new; it does not appear in Bradley's catalogue of formulae for Catalan's constant [4].

We have the approximation

$$G = 0.920558458... - \sum_{n=1}^{\infty} \frac{1}{n} \frac{1}{(16n^2 - 1)^2}$$

and since $G = 0.915965594...$ we note that

$$\sum_{n=1}^{\infty} \frac{1}{n} \frac{1}{(16n^2 - 1)^2} = 0.004592864...$$

We see that

$$G \approx 3(1 - \log 2) - \frac{1}{256} \sum_{n=1}^{\infty} \frac{1}{n^5} = 3(1 - \log 2) - \frac{\varsigma(5)}{256}$$

A better approximation is obtained by dealing separately with the first term of the series as follows

$$G \approx 3(1 - \log 2) - \frac{1}{225} - \frac{1}{256} \sum_{n=2}^{\infty} \frac{1}{n^5}$$

$$G \approx 3(1 - \log 2) - \frac{1}{225} - \frac{1}{256} [\varsigma(5) - 1]$$

$$\varsigma(5) = 1.0369277551...$$

It is interesting to note the appearance of $\varsigma(5)$ in these approximations.

*Mathematica* [8] also produces the result

(17) $$G = 1 - 16 \sum_{n=1}^{\infty} \frac{n}{(16n^2 - 1)^2}$$

In fact, this follows from the previous analysis because differentiating (5) gives us

$$\psi'(1+x) - \psi'(1-x) = -4x^3 \sum_{n=1}^{\infty} \frac{1}{n} \frac{1}{(n^2 - x^2)^2} - 4x \sum_{n=1}^{\infty} \frac{1}{n} \frac{1}{n^2 - x^2}$$

$$= -4x \sum_{n=1}^{\infty} \frac{n}{(n^2 - x^2)^2}$$

$$\psi'(1+x) - \psi'(1-x) = -4x \sum_{n=1}^{\infty} \frac{n}{(n^2 - x^2)^2}$$

A further differentiation results in

$$\psi''(1+x)+\psi''(1-x) = -8x^2 \sum_{n=1}^{\infty} \frac{n}{(n^2-x^2)^3} - 4\sum_{n=1}^{\infty} \frac{n}{(n^2-x^2)^2}$$

$$= -4\sum_{n=1}^{\infty} \frac{n^3}{(n^2-x^2)^3}$$

**Apéry's constant**

We may also use the above methodology to obtain a new series for $\zeta(3)$.

Differentiating (9) results in

(18)
$$\psi''\left(1+\frac{p}{2}\right)+\psi''\left(1-\frac{p}{2}\right) = -64p^4 \sum_{n=1}^{\infty} \frac{1}{n} \frac{1}{(4n^2-p^2)^3} - 80p^2 \sum_{n=1}^{\infty} \frac{1}{n} \frac{1}{(4n^2-p^2)^2}$$

$$-16\sum_{n=1}^{\infty} \frac{1}{n} \frac{1}{4n^2-p^2}$$

With $p=1/2$ in this we get

$$\psi''\left(1+\frac{1}{4}\right)+\psi''\left(\frac{3}{4}\right) = -64\frac{2^6}{2^4}\sum_{n=1}^{\infty} \frac{1}{n}\frac{1}{(16n^2-1)^3} - 80\frac{2^4}{2^2}\sum_{n=1}^{\infty} \frac{1}{n}\frac{1}{(16n^2-1)^2}$$

$$-16 \cdot 2^2 \sum_{n=1}^{\infty} \frac{1}{n} \frac{1}{16n^2-1}$$

We have using (10)

$$\psi''(1+x) - \psi''(x) = \frac{2}{x^3}$$

so that

$$\psi''\left(1+\frac{1}{4}\right) = \psi''\left(\frac{1}{4}\right) + 2^7$$

and we therefore obtain

$$\psi''\left(\frac{1}{4}\right)+\psi''\left(\frac{3}{4}\right)+2^7 = -64\frac{2^6}{2^4}\sum_{n=1}^{\infty}\frac{1}{n}\frac{1}{(16n^2-1)^3}-80\frac{2^4}{2^2}\sum_{n=1}^{\infty}\frac{1}{n}\frac{1}{(16n^2-1)^2}$$

$$-16\cdot 2^2\sum_{n=1}^{\infty}\frac{1}{n}\frac{1}{16n^2-1}$$

Hence we have

$$-112\varsigma(3)+2^7 = -64\frac{2^6}{2^4}\sum_{n=1}^{\infty}\frac{1}{n}\frac{1}{(16n^2-1)^3}-80\frac{2^4}{2^2}\sum_{n=1}^{\infty}\frac{1}{n}\frac{1}{(16n^2-1)^2}-16\cdot 2^2\sum_{n=1}^{\infty}\frac{1}{n}\frac{1}{16n^2-1}$$

which simplifies to

(19) $\quad 16\sum_{n=1}^{\infty}\frac{1}{n}\frac{1}{(16n^2-1)^3} = 7\varsigma(3)+48\log 2+20G-60$

where we used (7) and (16).

It may be noted that *Mathematica* directly confirms equations (16) and (19).

Higher derivatives of (18) will result in similar series for the Riemann zeta function $\varsigma(2n+1)$ and the Dirichlet beta function $\beta(2n)$.

The *Mathematica* evaluations for $\sum_{n=1}^{\infty}\frac{1}{n}\frac{1}{(16n^2-1)^r}$ for $r=4,5$ are set out below

```
Sum[1/(n (16 n^2-1)^4),{n,1,Infinity}]=

(-PolyGamma[3,1/4]+PolyGamma[3,3/4])/6144

+1/32 (140-44 Catalan-32 Log[8]-21 Zeta[3])

~0.0000197856927278423...

Sum[1/(n (16 n^2-1)^5),{n,1,Infinity}]=

(7 (PolyGamma[3,1/4]-PolyGamma[3,3/4]))/24576

+1/256 (-1260+372 Catalan+256 Log[8]+203 Zeta[3]

+31 Zeta[5])
~ 1.3173820678770678`^10-6
```

Rivoal and Zudilin [12] proved in 2003 that at least one of the seven numbers $\beta(2)$, $\beta(4)$, $\beta(6)$, $\beta(8)$, $\beta(10)$, $\beta(12)$ and $\beta(12)$ is irrational. I initially wondered whether it would be possible to use the series in (16) to show that $G - 3(1 - \log 2)$ is irrational. It is of course well known that $\log 2$ is irrational [9]. However, since this converges like $\varsigma(5)$, it seems that the rate of convergence will not be sufficiently rapid to demonstrate irrationality.

Apéry ([1] and [4]) showed in 1978 that $\varsigma(3)$ is irrational and hence we see that

$$4\sum_{n=1}^{\infty} \frac{1}{n} \frac{1}{(16n^2 - 1)^3} - 12\log 2 - 5G$$

is also irrational.

**ACKNOWLEDGEMENTS**

I thank R.J. Hughes for drawing my attention to identity (5) and also for verifying the *Mathematica* evaluations of (16), (17) and (19).

Donal F. Connon
Elmhurst
Dundle Road
Matfield
Kent TN12 7HD
dconnon@btopenworld.com